\documentclass{amsart}
\usepackage{amssymb}
\usepackage{amsmath}

\usepackage[pdftex,colorlinks,backref,pagebackref,hypertexnames=false,
            linkcolor=blue,urlcolor=blue,citecolor=blue]{hyperref}

\newcommand{\Z}{{\mathbb{Z}}}

\newcommand{\C}{{\mathbb{C}}}
\newcommand{\F}{{\mathbb{F}}}

\newcommand{\GL}{{\operatorname{GL}}}

\newcommand{\Irr}{\operatorname{Irr}}

\newcommand{\Stab}{\operatorname{Stab}}
\newcommand{\leg}{\operatorname{leg}}
\newcommand{\arm}{\operatorname{arm}}
\newcommand{\height}{\operatorname{ht}}

\newcommand{\dfour}[4]{\begin{footnotesize}\ensuremath{\hspace{-0.1cm}\begin{array}{ccc}&#2&\\#1\hspace{-0.2cm}&#3&\hspace{-0.2cm}#4\end{array}\hspace{-0.1cm}}\end{footnotesize}}

\renewcommand{\epsilon}{\varepsilon}

\newtheorem{theorem}{Theorem}[section]
\newtheorem{lemma}[theorem]{Lemma}
\newtheorem{proposition}[theorem]{Proposition}
\newtheorem{corollary}[theorem]{Corollary}

\newtheorem{defi}[theorem]{Definition}

\begin{document}

\title[Characters of the Sylow $p$-subgroups of $D_4(p^n)$]{Characters of the 
Sylow $\mathbf{p}$-subgroups of\\ the Chevalley groups $\mathbf{D_4(p^n)}$}

\author{Frank Himstedt, Tung Le and Kay Magaard}

\address{F.H.: Technische Universit\"at M\"unchen, Zentrum Mathematik --
         M11, Boltzmannstr. 3, 85748 Garching, Germany}
\address{T.L.: Mathematical Sciences, King's College, University of
  Aberdeen, Aberdeen AB24 3UE, U.K.}
\address{K.M.: School of Mathematics, University of Birmingham,
  Edgbaston, Birmingham B15 2TT, U.K.}

\email{F.H.: himstedt@ma.tum.de}
\email{T.L.: t.le@abdn.ac.uk}
\email{K.M.: k.magaard@bham.ac.uk}

\date{November 5, 2009}

% \keywords{Finite groups of Lie type, Sylow subgroups, Irreducible
% character degree}

% \subjclass[2000]{Primary 20C33, 20C15}

\begin{abstract}
Let $U(q)$ be a Sylow $p$-subgroup of the Chevalley groups $D_4(q)$
where $q$ is a power of a prime $p$. We describe a construction of all
complex irreducible characters of $U(q)$ and obtain a classification
of these irreducible characters via the root subgroups which are 
contained in the center of these characters. Furthermore, we show that
the multiplicities of the degrees of these irreducible characters are
given by polynomials in $q-1$ with nonnegative coefficients.
\end{abstract}

\maketitle

%%%%%%%%%%%%%%%%%%%%%%%%%%%%%%%%%%%%%%%%%%%%%%%%%%%%%%%%%%%%%%%%%%%%%%%%%%%%%%%%

 \section{Introduction}

Let $q$ be a power of a prime $p$ and $\F_q$ a field with
$q$ elements. The group $U_n(q)$ of all upper triangular $n \times
n$-matrices over $\F_q$ with all diagnonal entries equal to $1$ is a
Sylow $p$-subgroup of $\GL_n(\F_q)$. It was conjectured by
G.~Higman~\cite{Higman} that the number of conjugacy classes of
$U_n(q)$ is given by a polynomial in $q$ with integer coefficients.

Higman's conjecture was refined using the (complex) character theory
of $U_n(q)$. I.M.~Isaacs~\cite{IsaacsCharAlg} showed that the degrees
of the irreducible characters of $U_n(q)$ are of the form $q^e$, 
$0 \le e \le \mu(n)$ where the upper bound $\mu(n)$ is known
explicitly. G.~Lehrer~\cite{Lehrer} conjectured that the  
numbers $N_{n,e}(q)$ of irreducible characters of $U_n(q)$ of degree
$q^e$ are given by a polynomial in $q$ with integer
coefficients. I.M.~Isaacs suggested a strengthened form of
Lehrer's conjecture stating that $N_{n,e}(q)$ is given by a polynomial
in $q-1$ with nonnegative integer coefficients. So, Isaac's
conjecture implies Higman's and Lehrer's conjectures.

It is natural to consider these questions in the context of finite
groups of Lie type. Let $G(q)$ be a Chevalley group defined over
$\F_q$ and let $U(q)$ be a maximal unipotent subgroup of $G(q)$ such
that $U(q)$ is a Sylow $p$-subgroup of $G(q)$. S.~Goodwin and G.~R\"ohrle
\cite{GoodwinRoehrle} developed an algorithm for parametrizing the
conjugacy classes of $U(q)$ which is valid when $p$ is good for the
underlying root system. Using an implementation of this algorithm in
the computer algebra system GAP~\cite{GAP4}, they determined the
number of conjugacy classes of $U(q)$ for all Chevalley groups $G(q)$
of rank at most~$6$. As a consequence, they were able to show that
this number is given by a polynomial in $q$. 

A promising approach to an understanding of the complex irreducible
characters of $U(q)$ is the general concept of supercharacters which
was developed by P.~Diaconis and I.M.~Isaacs~\cite{DiaconisIsaacsSuper}.  
C.~Andr{\'e} and A.M.~Neto defined and studied supercharacters of
$U(q)$ if $G(q)$ if a classical group of type $B_n$, $C_n$ or $D_n$
and $p \ge 3$. These supercharacters are defined as products of
so-called elementary characters. Under the additional assumption 
$p \ge 2n$ %and the method of coadjoint orbits introduced by A.~Kirillov 
Andr{\'e} and Neto were able to show that every irreducible
character of $U(q)$ is a constituent of a unique supercharacter. In
particular, the supercharacters induce a partition on the set of
irreducible characters of $U(q)$. 

In this article we study the irreducible characters of a maximal
unipotent subgroup $U(q)$, that is a Sylow $p$-subgroup, of the
Chevalley groups $G(q)$ of type $D_4$. We describe a construction of
all irreducible characters of $U(q)$ which also works for bad
characteristic $p=2$. 
Our main result Theorem~\ref{thm:irrU} is a classification of the
irreducible characters of $U(q)$ via the root subgroups which are
contained in the center of the irreducible characters. 
In this way we get a natural partition of the set of irreducible
characters of $U(q)$ into families.
With each positive root we associate $q-1$ distinct irreducible
characters of $U(q)$ which we call \textit{midafis}. They are a
fundamental tool in the proof of our classification result and in some
sense they play a role similar to Andr{\'e}'s and Neto's elementary
characters.
We also obtain the degrees and for each degree the number of
irreducible characters of $U(q)$ of this degree. From this we can
conclude that an analogue of Isaac's conjecture is true for $U(q)$,
even in bad characteristic $p=2$. 

There are several differences between Andr{\'e}'s and Neto's methods
and ours:
\begin{itemize}
\item We do not make use of the natural matrix representation of the
  classical group of type $D_4$. Instead we use Lie theory and the
  underlying root system whenever it is possible.

\item Our midafis are always irreducible characters while Andr{\'e}'s
  and Neto's elementary characters are not necessarily irreducible.

\item We can describe a construction of all irreducible characters of
  $U(q)$. 

\item Our methods also work in bad characteristic $p=2$.
\end{itemize}

There are two reasons why we consider the Chevalley groups of type
$D_4$. First, they are not too far away from the classical
$A_n$-case. And second, the maximal unipotent subgroups of the
Chevalley groups of type $D_4$ are isomorphic to factor groups of the
maximal unipotent subgroups of the exceptional groups of
type $E_6$, $E_7$ and $E_8$.  
Therefore, we hope that this paper is a step towards the
classification of the irreducible characters of the maximal unipotent
subgroups of the ``large'' exceptional groups.

This paper is organized as follows: In Section~\ref{sec:nota}, we
introduce the general setup and fix notation. In
Section~\ref{sec:midafis}, we develop several tools which are essential
for our construction of the irreducible characters of $U(q)$: hook
subgroups and midafis. Finally, in Section~\ref{sec:IrrU}, we apply
these tools to obtain a construction and classification of the
irreducible characters of $U(q)$. Also in this section we determine
the degrees and numbers of the irreducible characters of $U(q)$ and
prove an analogue of Isaac's conjecture for $U(q)$.

%%%%%%%%%%%%%%%%%%%%%%%%%%%%%%%%%%%%%%%%%%%%%%%%%%%%%%%%%%%%%%%%%%%%%%%%%%%%%%%%

\section{Notation and Setup}
\label{sec:nota}

In this section, we introduce the setup and notation which will be
used throughout this paper.

\subsection{\texorpdfstring{Root system of type $\mathbf{D_4}$}{Root
system of type D4}} \label{subsec:d4} 
Let $\Phi$ be a root system of type $D_4$ in some Eucli\-dean space,
with basis $\Delta = \{\alpha_1, \alpha_2, \alpha_3, \alpha_4\}$ of
simple roots such that $\alpha_1, \alpha_2, \alpha_4$ are orthogonal
to each other. The Dynkin diagram of $\Phi$ is

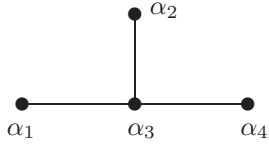
\begin{center}
\setlength{\unitlength}{1cm}
\begin{picture}(4,2)
\thinlines
\put(0.5,0.5){\circle*{0.17}}
\put(2,0.5){\circle*{0.17}}
\put(3.5,0.5){\circle*{0.17}}
\put(2,1.7){\circle*{0.17}}
\put(0.5,0.5){\line( 1, 0){1.5}}
\put(2,0.5){\line( 1, 0){1.5}}
\put(2,1.7){\line( 0, -1){1.2}}
\put(0.3,0.1){$\alpha_1$}
\put(1.9,0.1){$\alpha_3$}
\put(3.4,0.1){$\alpha_4$}
\put(2.2,1.7){$\alpha_2$}
\end{picture}
\end{center}

The positive roots are those roots which can be written as linear
combinations of the simple roots $\alpha_1, \alpha_2, \alpha_3,
\alpha_4$ with nonnegative coefficients and we write $\Phi_+$ for the
set of positive roots. We use the notation \dfour1121 for the 
root $\alpha_1 + \alpha_2 + 2 \alpha_3 + \alpha_4$ and we use a
similar notation for the remaining positive roots. The $12$ positive 
roots of $\Phi$ are given in Table~\ref{tab:posroots}. 

\begin{table}[!ht] 
\caption{Positive roots of the root system $\Phi$ of type $D_4$.} 
\label{tab:posroots}

\begin{center}
\begin{tabular}{c|llll}
\hline
\rule{0cm}{0.4cm}
Height & Roots &&&
\rule[-0.1cm]{0cm}{0.4cm}\\
\hline% \cline{1-3} \hline
\rule{0cm}{0.5cm}
5 & $\alpha_{12} := $ \dfour1121 &&&
\rule[-0.2cm]{0cm}{0.4cm}\\
\hline
\rule{0cm}{0.4cm}
4 & $\alpha_{11} := $ \dfour1111 &&&
\rule[-0.2cm]{0cm}{0.4cm}\\
\hline
\rule{0cm}{0.4cm}
3 & $\alpha_8 := $ \dfour1110 & $\alpha_9 := $ \dfour1011 & $\alpha_{10} := $ \dfour0111 &
\rule[-0.2cm]{0cm}{0.4cm}\\
\hline
\rule{0cm}{0.4cm}
2 & $\alpha_5 := $ \dfour1010 & $\alpha_6 := $ \dfour0110 & $\alpha_7 := $ \dfour0011 &
\rule[-0.2cm]{0cm}{0.4cm}\\
\hline
\rule{0cm}{0.4cm}
1 & $\alpha_1$  & $\alpha_2$  & $\alpha_3$   & $\alpha_4$ 
\rule[-0.2cm]{0cm}{0.4cm}\\
\hline
\end{tabular}
\end{center}
\end{table}

The numbering $\alpha_1, \alpha_2, \dots, \alpha_{12}$ of these 
roots is in accordance with the output of the CHEVIE \cite{CHEVIE} command
\begin{center}
\verb+CoxeterGroup("D", 4);+
\end{center}
We say that a nonempty subset $S \subseteq \Phi$ is \textit{closed} if
for all $\alpha, \beta \in S$ we have $\alpha + \beta \in S$ or 
$\alpha + \beta \not\in \Phi$. We write $\height(\alpha)$ for the height
of a positive root $\alpha$.

\subsection{\texorpdfstring{Chevalley groups of type
$\mathbf{D_4}$}{Chevalley groups of type D4}} \label{subsec:chevD4} 
Let $L$ be a simple complex Lie algebra with root system $\Phi$.
We~choose a Chevalley basis $\{h_r | r \in \Delta\} \cup \{e_r | r \in \Phi\}$ 
of $L$ auch that the structure constants $N_{rs}$ in 
$[e_r, e_s] = N_{rs} e_{r+s}$ are positive for all extraspecial pairs of roots 
$(r, s) \in \Phi \times \Phi$, see \cite[Section~4.2]{Carter1}. Fix a power
$q$ of some prime $p$ and let $G=G(q)=D_4(q)$ be the Chevalley group of type $D_4$
over the field $\F_q$ constructed from $L$, see \cite[Section~4.4]{Carter1}.
The group 
$$G = D_4(q) \cong P\Omega_8^+(q)$$ 
is simple and is generated by the root elements $x_\alpha(t)$ for
$\alpha \in \Phi$ and $t \in \F_q$. Let 
$X_\alpha := \langle x_\alpha(t) \, | \, t \in \F_q \rangle$ be the
root subgroup corresponding to $\alpha \in \Phi$. 
For positive roots, we use the abbreviation
$x_i(t) := x_{\alpha_i}(t)$, $i=1,2,\dots,12$. The 
commutators $[x_i(t), x_j(u)] = x_i(t)^{-1} x_j(u)^{-1} x_i(t) x_j(u)$ 
are given in Table~\ref{tab:commrelD4}. All 
$[x_i(t), x_j(u)]$ not listed in this table are equal to~$1$. 

\begin{table}[!ht] 
\caption{Commutator relations for type $D_4$.} 
\label{tab:commrelD4}

\begin{tabular}{llllll}
$\left[x_1(t), x_3(u)\right]$ & $=$ & $x_5(tu)$, \, &
$\left[x_1(t), x_6(u)\right]$ & $=$ & $x_8(-tu)$, \\
$\left[x_1(t), x_7(u)\right]$ & $=$ & $x_9(tu)$, \, &
$\left[x_1(t), x_{10}(u)\right]$ & $=$ & $x_{11}(-tu)$, \\
$\left[x_2(t), x_3(u)\right]$ & $=$ & $x_6(tu)$, \, &
$\left[x_2(t), x_5(u)\right]$ & $=$ & $x_8(-tu)$,\\
$\left[x_2(t), x_7(u)\right]$ & $=$ & $x_{10}(tu)$, \, &
$\left[x_2(t), x_9(u)\right]$ & $=$ & $x_{11}(-tu)$,\\
$\left[x_3(t), x_4(u)\right]$ & $=$ & $x_7(tu)$, \, &
$\left[x_3(t), x_{11}(u)\right]$ & $=$ & $x_{12}(-tu)$,\\
$\left[x_4(t), x_5(u)\right]$ & $=$ & $x_9(-tu)$, \, &
$\left[x_4(t), x_6(u)\right]$ & $=$ & $x_{10}(-tu)$,\\
$\left[x_4(t), x_8(u)\right]$ & $=$ & $x_{11}(-tu)$, \, &
$\left[x_5(t), x_{10}(u)\right]$ & $=$ & $x_{12}(-tu)$,\\
$\left[x_6(t), x_9(u)\right]$ & $=$ & $x_{12}(-tu)$,\, &
$\left[x_7(t), x_8(u)\right]$ & $=$ & $x_{12}(tu)$
\end{tabular}
\end{table}

Let $U=U(q)$ be the subgroup of $G=G(q)$ generated by the
elements $x_i(t)$ for $i=1,2, \dots, 12$ and $t \in \F_q$. So $U$ is a
maximal unipotent subgroup and a Sylow $p$-subgroup of $G$. In this
paper, we are interested in the complex irreducible characters of
$U$. Note that each element of $u \in U$ can be written uniquely as
\[
u = x_1(d_1) x_2(d_2) \cdots x_{12}(d_{12})
\]
where $d_1, \dots, d_{12} \in \F_q$. The multiplication of the
elements of $U$ is described by the commutator relations.
The center $Z(U) = X_{\alpha_{12}}$ is elementary abelian of order~$q$. 

\subsection{Characters, induction and restriction}
\label{subsec:charindres}

For any finite group $H$ let $\Irr(H)$ be the set of complex
irreducible characters~of~$H$ and let $(\cdot,\cdot)_H$ or
$(\cdot,\cdot)$ be the usual scalar product on the space of class
functions of $H$. Let $\mathbf{1}_H$ or $\mathbf{1}$ denote the
trivial character of $H$. If $\chi$ is a character of a subgroup~$H_1$
of $H$, then we write $\chi^H$ for the induced character, and if
$\chi$ is a character of $H$, we write $\chi|_{H_1}$ for the
restriction of $\chi$ to the subgroup $H_1$. 

\begin{defi}
Let $H$ be a finite group. We say that $\chi \in \Irr(H)$ is
\textit{almost faithful} if $Z(H) \not \subseteq \ker(\chi)$.
\end{defi}

Note that if $q$ is not prime, then the center $Z(U)$ is not
cyclic and in this case~$U$ does not have any faithful irreducible
characters. The almost faithful irreducible characters of $U$
are in some sense closest to being faithful.

For a field $K$ let $K^\times$ be its multiplicative group. 
In the whole paper, we fix a nontrivial linear character $\phi$ of
the group $(\F_q, +)$. So for $\alpha \in \Phi_+$ and $s \in \F_q$,
the map $\varphi_{\alpha, s}: X_\alpha \rightarrow \C^\times,
x_\alpha(d) \mapsto \phi(s \cdot d)$ is a linear character of 
the root subgroup~$X_\alpha$, and all irreducible characters of
$X_\alpha$ arise in this way.

%%%%%%%%%%%%%%%%%%%%%%%%%%%%%%%%%%%%%%%%%%%%%%%%%%%%%%%%%%%%%%%%%%%%%%%%%%%%%%%%

\section{Hook subgroups and midafis}
\label{sec:midafis}

For each positive root $\alpha$ we construct $q-1$ distinct
irreducible characters of $U$, called \textit{midafis}, which play a 
fundamental role in the classification and construction of the
irreducible characters of $U$ in this paper.

\subsection{Hook subgroups} \label{subsec:hooks}

With each positive root $\alpha$, we associate certain sets of
positive roots and a certain subgroup of $U$ which will be
used to define the midafis. 

\begin{defi} \label{def:hooks}
Let $\alpha \in \Phi_+$ be a positive root.
\begin{enumerate}
\item[(a)] The set $h_\alpha := \{\gamma \in \Phi_+ \, | \,
\text{there is   } \gamma' \in \Phi_+ \cup \{0\} \text{   such that   }
\gamma+\gamma'=\alpha\} \subseteq \Phi_+$ is called the \textit{hook}
corresponding to $\alpha$. 

\item[(b)] The subgroup $H_\alpha := \prod_{\gamma \in h_\alpha}
  X_\alpha \subseteq U$ is called the \textit{hook subgroup} of $U$
  corresponding to $\alpha$.

\item[(c)] We call
\[
\arm(h_\alpha) := \begin{cases}(h_\alpha \cap h_{\alpha_{12}})
  \setminus \{\alpha\} & \text{, if   } \alpha \neq \alpha_{12},\\
\{\alpha_{8}, \alpha_{9}, \alpha_{10}, \alpha_{11}\} & \text{, if   }
\alpha = \alpha_{12}\end{cases}
\]
the \textit{arm} and $\leg(h_\alpha) := h_\alpha \setminus
(\arm(h_\alpha) \cup \{\alpha\})$ the \textit{leg} of the hook $h_\alpha$.
\end{enumerate}
\end{defi}

\noindent \textbf{Remark:} (a) By the commutator relations, the hook
subgroup $H_\alpha$ does not depend on the order of the root subgroups
in the product.  

\smallskip

\noindent (b) Also, the commutator relations imply that the each hook
subgroup $H_\alpha$ is a special $p$-group of type $q^{1+2 \cdot {|\leg(h_\alpha)|}}$.

\medskip

\noindent \textbf{Example:} The hook subgroup
\[
H_{\alpha_{12}} = X_{\alpha_{3}} X_{\alpha_{5}} X_{\alpha_{6}}
X_{\alpha_{7}} X_{\alpha_{8}} X_{\alpha_{9}} 
X_{\alpha_{10}} X_{\alpha_{11}} X_{\alpha_{12}} 
\]
is the unipotent radical of the maximal parabolic subgroup of $G$
corresponding to the set $\{\alpha_1, \alpha_2, \alpha_4\}$ of simple
roots. So, $H_{\alpha_{12}}$ is a normal subgroup of $U$. Using the
notation from Subsection~\ref{subsec:d4}, we can picture
$H_{\alpha_{12}}$ as follows:

\begin{center}
\begin{tabular}{ccccc}
\dfour1110 & \dfour1011 & \dfour 0111 & \dfour1111 & \dfour1121 

\medskip\\

&&&& \dfour0011

\smallskip\\

&&&& \dfour0110

\medskip\\

&&&& \dfour1010

\medskip\\

&&&& \dfour0010
\end{tabular}
\end{center}

\subsection{Midafis} \label{subsec:midafis}

For positive roots $\alpha \in \Phi_+$ we set
$V_\alpha := \prod_{\gamma \in \Phi_+ \setminus \leg(\alpha)} X_\gamma$.
The midafis of $U$ associated with $\alpha$ will be characters which
are induced from $V_\alpha$. We are going to use the following lemma:

\begin{lemma} \label{la:basegrp}
Let $\alpha \in \Phi_+$ and $s \in \F_q\times$.
\begin{enumerate}
\item[(a)] The subset $\Phi_+ \setminus \leg(\alpha) \subseteq \Phi$
  is closed. 

\item[(b)] The set $V_\alpha$ is a subgroup of $U$.

\item[(c)] We have $X_\alpha \cap [V_\alpha, V_\alpha] = \{1\}$.

\item[(d)] For each $s \in \F_q^\times$ there is a linear 
  $\lambda_{\alpha,s} \in \Irr(V_\alpha)$ such that
  $\lambda_{\alpha,s}|_{X_\alpha} = \varphi_{\alpha,s}$ and
  $X_\gamma \subseteq \ker(\lambda_{\alpha,s})$ for all 
  $\gamma \in \Phi_+$ with $\height(\gamma) > \height(\alpha)$.
\end{enumerate}
\end{lemma}
\begin{proof}
\noindent (a) Suppose not. So there is $\beta \in \leg(\alpha)$ and  
$\gamma, \gamma' \in \Phi_+ \setminus \leg(\alpha)$ such that  
$\beta = \gamma+\gamma'$. But then $\gamma \in \leg(\alpha)$ or
$\gamma' \in \leg(\alpha)$, a contradiction.

\smallskip

\noindent (b) follows from (a).

\smallskip

\noindent (c) Suppose 
$X_\alpha \cap [V_\alpha, V_\alpha] \neq \{1\}$. The commutator
relations imply that there are 
$\gamma, \gamma' \in \Phi_+ \setminus \leg(\alpha)$ such that
$\gamma+\gamma'=\alpha$. But then $\gamma \in \leg(\alpha)$ or
$\gamma' \in \leg(\alpha)$, which is a contradiction.

\smallskip

\noindent (d) follows from (c).
\end{proof}

\begin{defi} \label{def:midafi}
Let $\alpha$ be a positive root. We call the characters
$\mu_{\alpha,s} := \lambda_{\alpha,s}^U$ for $s \in \F_q^\times$ the
\textit{midafis} of~$U$ associated with $\alpha$.
\end{defi}

\begin{proposition} \label{prop:midafi}
Let $\alpha \in \Phi_+$ be a positive root. The midafis 
$\mu_{\alpha,s}$ for $s \in \F_q^\times$ are $q-1$ distinct
irreducible characters of $U$ and 
$\mu_{\alpha,s}|_{X_\alpha} = \mu_{\alpha,s}(1) \cdot \varphi_{\alpha,s}$.
\end{proposition}
\begin{proof}
Note that for all $\gamma \in \Phi_+$ with $\height(\gamma) >
\height(\alpha)$, we have $X_\gamma \subseteq \ker(\lambda_{\alpha,s})$.
So the statement about the restriction $\mu_{\alpha,s}|_{X_\alpha}$ is clear
by the commutator relations and the definition of induced characters,
and from this we immediately get $\mu_{\alpha,s} \neq \mu_{\alpha,s'}$
for $s \neq s' \in \F_q^\times$. So we only have to show that
$\mu_{\alpha,s}$ is irreducible.  

Let 
$\leg(\alpha) = \{\alpha_{i_1}, \alpha_{i_2}, \dots, \alpha_{i_m}\}$
with $\height(\alpha_{i_1}) \ge \height(\alpha_{i_2}) \ge \dots
\ge \height(\alpha_{i_m})$. From the commutator relations we get
\[
V_\alpha \unlhd X_{\alpha_1} V_\alpha \unlhd X_{\alpha_1} X_{\alpha_2}
V_\alpha \unlhd \dots \unlhd \left(\prod_{j=1}^m X_{\alpha_j}\right) V_\alpha=U.
\]
Considering character values on $X_{\alpha-\alpha_{i_1}}$ we see
that the inertia subgroup of $\lambda_{\alpha,s}$ in $X_{\alpha_1}
V_\alpha$ is equal to $V_\alpha$. Hence by Clifford
theory~\cite[Theorem~(6.11)]{Isaacs:76}, the induced character 
$\lambda_{\alpha,s}^{X_{\alpha_{i_1}} V_\alpha}$ is irreducible.
Similarly, considering character values on $X_{\alpha-\alpha_{i_2}}$
we see that the inertia subgroup of 
$\lambda_{\alpha,s}^{X_{\alpha_{i_1}} V_\alpha}$ in 
$X_{\alpha_{i_1}} X_{\alpha_{i_2}} V_\alpha$ is equal to 
$X_{\alpha_{i_1}} V_\alpha$. Again, \cite[Theorem~(6.11)]{Isaacs:76}
implies the irreducibility of the induced character 
$\lambda_{\alpha,s}^{X_{\alpha_{i_1}} X_{\alpha_{i_2}} V_\alpha}$.
Continuing in this way, we get that
$\mu_{\alpha,s}=\lambda_{\alpha,s}^U$ is irreducible. 
\end{proof}

\subsection{Hook subgroups for the group of upper triangular matrices} 
\label{subsec:hooksGLn}

The notation of hooks, arms, legs and midafis is motivated by
root systems of type $A$, that is, by the structure of the Sylow
$p$-subgroups of $\GL_n(\F_q)$. The group $U_n(q)$ of all upper
unitriangular matrices over $\F_q$ is a Sylow $p$-subgroup of
$\GL_n(\F_q)$. The root system of $\GL_n(\F_q)$ with respect to the
maximal torus of diagonal matrices has simple roots 
$\alpha_1, \alpha_2, \dots, \alpha_{n-1}$ such that the nodes 
corresponding to $\alpha_i$ and $\alpha_{i+1}$ are joined in the
Dynkin diagram for $i=1,2,\dots,n-2$. The positive roots are 
the roots $\alpha_{ij} := \alpha_i + \dots + \alpha_j$ 
for all $1 \le i \le j \le n-1$. The root subgroup $X_{\alpha_{ij}}$
consists of the matrices $I_n+ t \cdot e_{i,j+1}$ for $t \in \F_q$,
where $I_n$ is the $n \times n$-identity matrix and $e_{i,j+1}$ is the
$n \times n$-matrix with zero entries except a single entry $1$ in
position $(i, j+1)$. Hooks and hook subgroups can be defined for
$U_n(q)$ in the same way as in Definition~\ref{def:hooks} (a), (b). 
For $n=8$, the hooks subgroup corresponding to $\alpha_{2,6}$
can be pictured as follows:
\begin{small}
\[
\left(\begin{array}{cccccccc}
1 & . & . & . & . & . & . & .\\
\cline{3-6}
. & 1 & \multicolumn{1}{|c}{*} & * & * & \multicolumn{1}{|c|}{*} & . & .\\
\cline{3-6}
. & . & 1 & . & . & \multicolumn{1}{|c|}{*} & . & .\\
. & . & . & 1 & . & \multicolumn{1}{|c|}{*} & . & .\\
. & . & . & . & 1 & \multicolumn{1}{|c|}{*} & . & .\\
\cline{6-6}
. & . & . & . & . & 1 & . & .\\
. & . & . & . & . & . & 1 & .\\
. & . & . & . & . & . & . & 1
\end{array}\right)
\]
\end{small}
Associated with each positive root $\alpha_{ij}$ are $q-1$ irreducible
characters of $U_n(q)$ which can be constructed in a way analogous to
our midafis for type $D_4$. The word \textit{midafi} is an
abbreviation for \textit{minimal degree almost faithful irreducible}
which comes from the fact that the midafis of $U_n(q)$ can be
interpreted as almost faithful irreducible characters of minimal
degree of suitable factor groups of $U_n(q)$. For details,
see~\cite[Section~2.2]{TungPhD}.

%%%%%%%%%%%%%%%%%%%%%%%%%%%%%%%%%%%%%%%%%%%%%%%%%%%%%%%%%%%%%%%%%%%%%%%%%%%%%%%%

\section{\texorpdfstring{Irreducible characters of $U$}{Irreducible characters of U}}
\label{sec:IrrU}

In this section, we describe a construction of all irreducible
characters of $U$ for all prime powers $q$. In particular, we obtain
the number and the degrees of all irreducible characters of $U$. The
main result is the following theorem. 

\begin{theorem} \label{thm:irrU}
For every prime power $q$ the irreducible characters of $U$ are given
by Table~\ref{tab:irrU}. 
\end{theorem}

We begin with some comments on Table~\ref{tab:irrU}. There are
17 families of irreducible characters of $U$ and each row of Table~\ref{tab:irrU}
represents one of these families. The~first column gives notation for
these families of characters. Note that the family
$\mathcal{F}_{8,9,10}^{odd}$ exists only for odd $q$, while
$\mathcal{F}_{8,9,10}^{even}$ exists only if $q$ is even.
The second column of Table~\ref{tab:irrU} gives notation for~the
irreducible characters in each family. The first one, two or three
indices of this notation describe the positive roots $\alpha_j$ of
maximal height such that~$X_{\alpha_j}$ is not contained in the kernel
of the irreducible characters in the family. If there are two types of
characters in the family, there is an additional index which describes 
the degree of the characters. The remaining indices are
parameters which can take values from the parameter set in the third
column. The fourth column lists the number of irreducible characters
in the family and the last column gives their degrees. Note that we
use a slightly different notation for the family of linear characters.

\medskip

\noindent \textbf{Example:} Family $\mathcal{F}_8$ constists of 2 types of
characters. The characters $\chi_{8,q^3,a_1,a_2}$, where $a_1, a_2$
vary over $\F_q^\times$, are $(q-1)^2$ distinct irreducible characters
of degree $q^3$, and the characters $\chi_{8,q^2,a,b_1,b_2}$, where
$a$ varies over $\F_q^\times$ and $b_1,b_2$ vary over $\F_q$, are
$q^2(q-1)$ distinct irreducible characters of $U$ of degree~$q^2$.
Furthermore, one has 
$X_{\alpha_{9}} X_{\alpha_{10}} X_{\alpha_{11}}
X_{\alpha_{12}} \subseteq \ker(\chi_{8,q^3,a_1,a_2}),
\ker(\chi_{8,q^2,a,b_1,b_2})$ and 
$X_{\alpha_8} \not\subseteq \ker(\chi_{8,q^3,a_1,a_2})$ and 
$X_{\alpha_8} \not\subseteq \ker(\chi_{8,q^2,a,b_1,b_2})$. 

\begin{table}[!ht] 
\caption{Irreducible characters of $U$.} 
\label{tab:irrU}

\begin{center}
\begin{tabular}{l|l|l|l|l}
\hline
\rule{0cm}{0.4cm}
Family & Notation & Parameter set & Number & Degree
\rule[-0.1cm]{0cm}{0.4cm}\\
\hline
\rule{0cm}{0.5cm}
$\mathcal{F}_{12}$ & $\chi_{12,a,b_1,b_2,b_3}$ & $\F_q^\times \times \F_q \times \F_q
\times \F_q$ & $q^3(q-1)$ & $q^4$
\rule[-0.2cm]{0cm}{0.4cm}\\
\hline
\rule{0cm}{0.4cm}
$\mathcal{F}_{11}$ & $\chi_{11,a,b_1,b_2,b_3,b_4}$ & $\F_q^\times \times \F_q \times \F_q
\times \F_q \times \F_q$ & $q^4(q-1)$ & $q^3$
\rule[-0.2cm]{0cm}{0.4cm}\\
\hline
\rule{0cm}{0.4cm}
$\mathcal{F}_{8,9,10}^{odd}$  & $\chi_{8,9,10,a_1,a_2,a_3,b}$ & $\F_q^\times
\times \F_q^\times \times \F_q^\times \times \F_q$ & $q(q-1)^3$ & $q^3$
\rule[-0.2cm]{0cm}{0.4cm}\\
\hline
\rule{0cm}{0.4cm}
$\mathcal{F}_{8,9,10}^{even}$  & $\chi_{8,9,10,q^3,a_1,a_2,a_3}$ & 
$\F_q^\times \times \F_q^\times \times \F_q^\times$ & $(q-1)^3$ & $q^3$\\
& $\chi_{8,9,10,\frac{q^3}{2},x,a_1,a_2,a_3,a_4}$ & $\Z_4 \times
\F_q^\times \times \F_q^\times \times \F_q^\times \times \F_q^\times$ &
$4(q-1)^4$ & $\frac{q^3}{2}$
\rule[-0.2cm]{0cm}{0.4cm}\\
\hline
\rule{0cm}{0.4cm}
$\mathcal{F}_{8,9}$ & $\chi_{8,9,q^3,a_1,a_2,a_3}$ & $\F_q^\times
\times \F_q^\times \times \F_q^\times$ & $(q-1)^3$ & $q^3$\\
& $\chi_{8,9,q^2,a_1,a_2,b_1,b_2}$ & $\F_q^\times
\times \F_q^\times \times \F_q \times \F_q$ & $q^2(q-1)^2$ & $q^2$
\rule[-0.2cm]{0cm}{0.4cm}\\
\hline
\rule{0cm}{0.4cm}
$\mathcal{F}_{8,10}$ & $\chi_{8,10,q^3,a_1,a_2,a_3}$ & $\F_q^\times
\times \F_q^\times \times \F_q^\times$ & $(q-1)^3$ & $q^3$\\
& $\chi_{8,10,q^2,a_1,a_2,b_1,b_2}$ & $\F_q^\times
\times \F_q^\times \times \F_q \times \F_q$ & $q^2(q-1)^2$ & $q^2$
\rule[-0.2cm]{0cm}{0.4cm}\\
\hline
\rule{0cm}{0.4cm}
$\mathcal{F}_{9,10}$ & $\chi_{9,10,q^3,a_1,a_2,a_3}$ & $\F_q^\times
\times \F_q^\times \times \F_q^\times$ & $(q-1)^3$ & $q^3$\\
& $\chi_{9,10,q^2,a_1,a_2,b_1,b_2}$ & $\F_q^\times
\times \F_q^\times \times \F_q \times \F_q$ & $q^2(q-1)^2$ & $q^2$
\rule[-0.2cm]{0cm}{0.4cm}\\
\hline
\rule{0cm}{0.4cm}
$\mathcal{F}_{8}$ & $\chi_{8,q^3,a_1,a_2}$ & $\F_q^\times
\times \F_q^\times$ & $(q-1)^2$ & $q^3$\\
& $\chi_{8,q^2,a_1,b_1,b_2}$ & $\F_q^\times
\times \F_q \times \F_q$ & $q^2(q-1)$ & $q^2$
\rule[-0.2cm]{0cm}{0.4cm}\\
\hline
\rule{0cm}{0.4cm}
$\mathcal{F}_{9}$ & $\chi_{9,q^3,a_1,a_2}$ & $\F_q^\times
\times \F_q^\times$ & $(q-1)^2$ & $q^3$\\
& $\chi_{9,q^2,a_1,b_1,b_2}$ & $\F_q^\times
\times \F_q \times \F_q$ & $q^2(q-1)$ & $q^2$
\rule[-0.2cm]{0cm}{0.4cm}\\
\hline
\rule{0cm}{0.4cm}
$\mathcal{F}_{10}$ & $\chi_{10,q^3,a_1,a_2}$ & $\F_q^\times
\times \F_q^\times$ & $(q-1)^2$ & $q^3$\\
& $\chi_{10,q^2,a_1,b_1,b_2}$ & $\F_q^\times
\times \F_q \times \F_q$ & $q^2(q-1)$ & $q^2$
\rule[-0.2cm]{0cm}{0.4cm}\\
\hline
\rule{0cm}{0.4cm}
$\mathcal{F}_{5,6,7}$ & $\chi_{5,6,7,a_1,a_2,a_3,b_1,b_2}$ &
$\F_q^\times \times \F_q^\times \times \F_q^\times \times \F_q \times
\F_q$ & $q^2(q-1)^3$ & $q$
\rule[-0.2cm]{0cm}{0.4cm}\\
\hline
\rule{0cm}{0.4cm}
$\mathcal{F}_{5,6}$ & $\chi_{5,6,a_1,a_2,b_1,b_2}$ &
$\F_q^\times \times \F_q^\times \times \F_q \times \F_q$ & $q^2(q-1)^2$ & $q$
\rule[-0.2cm]{0cm}{0.4cm}\\
\hline
\rule{0cm}{0.4cm}
$\mathcal{F}_{5,7}$ & $\chi_{5,7,a_1,a_2,b_1,b_2}$ &
$\F_q^\times \times \F_q^\times \times \F_q \times \F_q$ & $q^2(q-1)^2$ & $q$
\rule[-0.2cm]{0cm}{0.4cm}\\
\hline
\rule{0cm}{0.4cm}
$\mathcal{F}_{6,7}$ & $\chi_{6,7,a_1,a_2,b_1,b_2}$ &
$\F_q^\times \times \F_q^\times \times \F_q \times \F_q$ & $q^2(q-1)^2$ & $q$
\rule[-0.2cm]{0cm}{0.4cm}\\
\hline
\rule{0cm}{0.4cm}
$\mathcal{F}_{5}$ & $\chi_{5,a,b_1,b_2}$ &
$\F_q^\times \times \F_q \times \F_q$ & $q^2(q-1)$ & $q$
\rule[-0.2cm]{0cm}{0.4cm}\\
\hline
\rule{0cm}{0.4cm}
$\mathcal{F}_{6}$ & $\chi_{6,a,b_1,b_2}$ &
$\F_q^\times \times \F_q \times \F_q$ & $q^2(q-1)$ & $q$
\rule[-0.2cm]{0cm}{0.4cm}\\
\hline
\rule{0cm}{0.4cm}
$\mathcal{F}_{7}$ & $\chi_{7,a,b_1,b_2}$ &
$\F_q^\times \times \F_q \times \F_q$ & $q^2(q-1)$ & $q$
\rule[-0.2cm]{0cm}{0.4cm}\\
\hline
\rule{0cm}{0.4cm}
$\mathcal{F}_{lin}$ & $\chi_{lin,b_1,b_2,b_3,b_4}$ & $\F_q \times \F_q \times \F_q
\times \F_q$ & $q^4$ & $1$
\rule[-0.2cm]{0cm}{0.4cm}\\
\hline
\end{tabular}
\end{center}
\end{table}

\subsection{Proof of Theorem~\ref{thm:irrU}} \label{subsec:proofirrU}

We describe the definition and construction of the irreducible
characters in each family. Our main tools will be the midafis,
Clifford theory and the commutator relations.

\medskip

\noindent \textbf{The irreducible characters in family $\mathcal{F}_{12}$.} 
In this subsection we describe the construction of the irreducible
characters in family $\mathcal{F}_{12}$ which is the family of the
almost faithful irreducible characters of $U$. 

We have already seen in
Subsection~\ref{subsec:hooks} that the hook subgroup $H_{\alpha_{12}}$
is normal in $U$. Since $H_{\alpha_{12}}$ is a special $p$-group of
type $q^{1+8}$, it has has $q^8$ linear characters and $q-1$
irreducible characters of degree~$q^4$. 
We have $X_{\alpha_{12}} \not \subseteq \ker(\mu_{\alpha_{12},s})$ and 
$\mu_{\alpha_{12},s}|_{X_{\alpha_{12}}} \neq
\mu_{\alpha_{12},s'}|_{X_{\alpha_{12}}}$ for all 
$s \neq s' \in \F_q^\times$. So the restrictions of the midafis
$\mu_{\alpha_{12},s}$ are exactly the nonlinear irreducible characters
of~$H_{\alpha_{12}}$. In particular, each nonlinear
irreducible character of~$H_{\alpha_{12}}$ extends to $U$.

By Clifford theory, the almost faithful irreducible characters of $U$
are those lying over the nonlinear irreducible characters of $H_{\alpha_{12}}$.
The factor group $U/H_{\alpha_{12}}$ is isomorphic to 
$X_{\alpha_1} \times X_{\alpha_2} \times X_{\alpha_4} \cong \F_q
\times \F_q \times \F_q$. So by Gallagher's
theorem~\cite[Corollary~(6.17)]{Isaacs:76}, there are $q^3(q-1)$
almost faithful irreducible characters of $U$. They all have degree
$q^4$ and can be parametrized by $a \in \F_q^\times$ (parametrizing the
midafi) and $b_1, b_2, b_3 \in \F_q$ (parametrizing the irreducible
characters of $U/H_{\alpha_{12}}$). This proves all statements about
the family $\mathcal{F}_{12}$ in Table~\ref{tab:irrU}. 

\medskip

\noindent \textbf{The irreducible characters in family $\mathcal{F}_{11}$.} 
The characters in family $\mathcal{F}_{11}$ are the irreducible 
characters $\chi$ of $U$ such that $X_{\alpha_{12}} \subseteq \ker(\chi)$ 
and $X_{\alpha_{11}} \not \subseteq \ker(\chi)$.
We are going to work in the factor group 
$\overline{U} := U/Z(U) = U/X_{\alpha_{12}}$. By the commutator
relations, we have the following subnormal series:
\[
\{1\} \unlhd \overline{H}_{11} := H_{\alpha_{11}}X_{\alpha_{12}}/X_{\alpha_{12}} 
\unlhd \overline{N}_{11} := \prod_{\genfrac{}{}{0pt}{}{i=1}{i \neq 3}}^{12} X_{\alpha_i}/X_{\alpha_{12}} 
\unlhd \overline{U},
\]
where $\overline{N}_{11}$ does not depend on the order of the root
subgroups in the product. Furthermore, 
$\overline{H}_{11} \cong H_{\alpha_{11}}$ is a special $p$-group of
type $q^{1+6}$, the factor group 
$\overline{N}_{11}/\overline{H}_{11} \cong X_{\alpha_5} \times
X_{\alpha_6} \times X_{\alpha_7}$ is elementary abelian of order $q^3$
and $\overline{U} / \overline{N}_{11} \cong X_{\alpha_3}$ is
elementary abelian of order $q$.

First, we construct all $\psi \in \Irr(\overline{N}_{11})$ 
with $X_{\alpha_{11}} X_{\alpha_{12}} / X_{\alpha_{12}} \not \subseteq \ker(\psi)$.
Since~$\overline{H}_{11}$ is special of type $q^{1+6}$, it has $q^6$
linear characters and $q-1$ almost faithful irreducible characters of
degree $q^3$. For each $s \in \F_q^\times$, one has
$X_{\alpha_{12}} \subseteq \ker(\mu_{\alpha_{11},s})$. Hence we can
identify the midafi $\mu_{\alpha_{11},s}$ with an irreducible character of
$\overline{U}$, which we also denote by~$\mu_{\alpha_{11},s}$.
Since $\deg(\mu_{\alpha_{11},s})=q^3$ and 
$X_{\alpha_{11}} \not \subseteq \ker(\mu_{\alpha_{11},s})$, the
restriction of $\mu_{\alpha_{11},s}$ to $\overline{H}_{11}$ is
irreducible. Thus, also $\mu_{\alpha_{11},s}|_{\overline{N}_{11}}$ is
irreducible and each almost faithful irreducible character of
$\overline{H}_{11}$ extends to $\overline{N}_{11}$. Gallagher's
theorem~\cite[Corollary~(6.17)]{Isaacs:76} applied to 
$\overline{H}_{11} \unlhd \overline{N}_{11}$ implies that each almost
faithful irreducible character of $\overline{H}_{11}$ extends to
$\overline{N}_{11}$ in $q^3$ ways. Hence there are $q^3(q-1)$
irreducible characters of $\overline{N}_{11}$ such that 
$X_{\alpha_{11}} X_{\alpha_{12}} / X_{\alpha_{12}}$ is not contained
in their kernel, these characters have degree $q^3$ and can be
parametrized by $\psi_{a,b_1,b_2,b_3}$ where $a \in \F_q^\times$ and
$b_i \in \F_q$. 

Next, we show that each $\psi_{a,b_1,b_2,b_3}$ extends to
$\overline{U}$. Since 
$
X_{\alpha_3} \cdot \prod_{i \ge 5} X_{\alpha_i} / X_{\alpha_{12}}
$
is a normal abelian subgroup of $\overline{U}$ of index $q^3$, 
we can conclude from Ito's theorem \cite[Theorem (6.15)]{Isaacs:76}
that the degrees of all irreducible characters of $\overline{U}$
divide $q^3$. So from Clifford theory and Gallagher's theorem applied
to $\overline{N}_{11} \unlhd \overline{U}$ it follows that each 
$\psi_{a,b_1,b_2,b_3}$ extends to $\overline{U}$ in $q$ ways. We
denote these extensions by $\psi_{a,b_1,b_2,b_3,b_4}$. Inflating 
$\psi_{a,b_1,b_2,b_3,b_4}$ to $U$, we get the $q^4(q-1)$
irreducible characters $\chi_{a,b_1,b_2,b_3,b_4}$ of degree $q^3$ of~$U$. 

\medskip

\noindent \textbf{The irreducible characters in family
  $\mathcal{F}_{8,9,10}^{odd}$.} 
We assume that $q$ is odd. The characters in
$\mathcal{F}_{8,9,10}^{odd}$ are the irreducible characters $\chi$ of
$U$ such that $X_{\alpha_{i}} \subseteq \ker(\chi)$ for $i=11,12$
and $X_{\alpha_{j}} \not \subseteq \ker(\chi)$ for $j=8,9,10$.
We are going to work in the factor group 
$\overline{U} := U/X_{\alpha_{11}}X_{\alpha_{12}}$. Note that
$\overline{U}$ is a semidirect product 
$\overline{U} = \overline{K} \ltimes \overline{A}$
of the elementary abelian groups 
\[
\overline{K} = X_{\alpha_1} X_{\alpha_2} X_{\alpha_4} \quad \text{and}
\quad \overline{A} := X_{\alpha_3} \cdot \prod_{i \ge 5} X_{\alpha_i}
/ X_{\alpha_{11}}X_{\alpha_{12}}.
\]
We consider the action of $X_{\alpha_1} X_{\alpha_2} X_{\alpha_4}$ on 
$\Irr(\overline{A})$ by conjugation. For field elements 
$x,a,b,c,d,e,f \in \F_q$ we define a linear character 
$\lambda_{x,a,b,c,d,e,f} \in \Irr(\overline{A})$ by
\[
x_3(d_3) x_5(d_5) x_6(d_6) \cdots x_{10}(d_{10})
X_{\alpha_{11}} X_{\alpha_{12}} 
\mapsto \phi(x \cdot d_3 + a \cdot d_5 + b \cdot d_6 + \dots + f \cdot d_{10}),
\]
where $\phi$ is a nontrivial linear character of $(\F_q,+)$ as in
Subsection~\ref{subsec:charindres}. We claim that
\begin{equation} \label{eq:repset}
\{\lambda_{x,0,0,0,d,e,f} \, | \, x \in \F_q, d,e,f \in
\F_q^\times\}
\end{equation}
is a set of representatives for the action of 
$X_{\alpha_1} X_{\alpha_2} X_{\alpha_4}$ on the set $\lambda \in
\Irr(\overline{A})$ such that $\lambda|_{X_{\alpha_i}}$ is nontrivial
for $i=8,9,10$. From the commutator relations we obtain
\begin{eqnarray} \label{eq:comrel8910}
{^{x_1(r) x_2(s) x_4(t)}(x_3(d_3) x_5(d_5) \cdots x_{10}(d_{10})}) & =
& x_3(d_3) x_5(d_5+rd_3) \nonumber\\
&& \cdot x_6(d_6+sd_3)x_7(d_7-td_3) \nonumber\\
&& \cdot x_8(d_8-sd_5-rd_6-rsd_3)\\
&& \cdot x_9(d_9-td_5+rd_7-rtd_3) \nonumber\\
&& \cdot x_{10}(d_{10}-td_6+sd_7-std_3). \nonumber
\end{eqnarray}
So, $(\lambda_{x,0,0,0,d,e,f})^{x_1(r) x_2(s) x_4(t)}=\lambda_{x',0,0,0,d',e',f'}$ if and only if
\begin{eqnarray*}
\phi(x'd_3+d'd_8+e'd_9+f'd_{10}) & = & \phi((x-drs-ert-fst)d_3\\
&& +(-ds-et)d_5+(-dr-ft)d_6\\
&& +(er+fs)d_7+dd_8+ed_9+fd_{10})
\end{eqnarray*}
for all $d_i \in \F_q$. Thus, if 
$(\lambda_{x,0,0,0,d,e,f})^{x_1(r) x_2(s) x_4(t)}=\lambda_{x',0,0,0,d',e',f'}$
then the coefficients of $d_5, d_6, \dots, d_{10}$ imply $d'=d, e'=e,
f'=f$ and 
\[
\left(\begin{array}{ccc}
0 & d & e\\
d & 0 & f\\
e & f & 0
\end{array}\right) \cdot 
\left(\begin{array}{c}
r\\
s\\
t
\end{array}\right) = 
\left(\begin{array}{c}
0\\
0\\
0
\end{array}\right).
\]
The determinant of the coefficient matrix of this system of linear
equations is $2 def$. Since $q$ is odd, we conclude $r=s=t=0$ and so
$x'=x$. So, (\ref{eq:repset}) is indeed a set of representatives and 
$\Stab_{\overline{K}}(\lambda_{x,0,0,0,d,e,f}) = \{1\}$.
By Clifford theory, the induced characters 
$\psi_{a_1,a_2,a_3,b} := \lambda_{b,0,0,0,a_1,a_2,a_3}^{\overline{U}}$
are $q(q-1)^3$ distinct irreducible characters of
$\overline{U}$. By inflation, we obtain 
$\chi_{8,9,10,a_1,a_2,a_3,b} \in \Irr(U)$ of degree $q^3$ of~$U$.

\medskip

\noindent \textbf{The irreducible characters in family
  $\mathcal{F}_{8,9,10}^{even}$.} 
We assume that $q$ is even. The characters in
$\mathcal{F}_{8,9,10}^{even}$ are the irreducible characters $\chi$ of
$U$ such that $X_{\alpha_{i}} \subseteq \ker(\chi)$ for $i=11,12$
and $X_{\alpha_{j}} \not \subseteq \ker(\chi)$ for $j=8,9,10$.
Again, we are going to work in the factor group 
$\overline{U} := U/X_{\alpha_{11}}X_{\alpha_{12}}$. We define
subgroups $\overline{K}= X_{\alpha_1} X_{\alpha_2} X_{\alpha_4}$ and
$\overline{A}$ of $\overline{U}$ in the same way as for odd $q$ 
and obtain the decomposition $\overline{U} = \overline{K} \ltimes \overline{A}$. 
We consider the conjugation action of $X_{\alpha_1}X_{\alpha_2}X_{\alpha_4}$ 
on the set of irreducible characters of the elementary abelian
normal subgroup $\overline{A}$. As for odd $q$, we define linear 
characters $\lambda_{x,a,b,c,d,e,f} \in \Irr(\overline{A})$ by
\[
x_3(d_3) x_5(d_5) x_6(d_6) \cdots x_{10}(d_{10})
X_{\alpha_{11}} X_{\alpha_{12}} 
\mapsto \phi(x \cdot d_3 + a \cdot d_5 + b \cdot d_6 + \dots + f \cdot d_{10}),
\]
where $\phi$ is a nontrivial linear character of $(\F_q,+)$ and 
$x,a,b,c,d,e,f \in \F_q$. If $cdef \neq 0$, then 
$\{\, defz^2 + cdz | \, z \in \F_q\}$ is a subgroup of $(\F_q,+)$ of
index~$2$. Choose $t_{c,d,e,f} \in \F_q \setminus \{\, defz^2 + cdz |
\, z \in \F_q\}$. We claim that
\begin{equation} \label{eq:repseteven}
\{\lambda_{0,0,0,0,d,e,f} \, | \, d,e,f \in \F_q^\times\} \cup
\{\lambda_{x,0,0,c,d,e,f} \, | \, x \in \{0,t_{c,d,e,f}\}, c,d,e,f \in
  \F_q^\times\}
\end{equation}
is a set of representatives for the action of 
$X_{\alpha_1} X_{\alpha_2} X_{\alpha_4}$ on the set $\lambda \in
\Irr(\overline{A})$ such that $\lambda|_{X_{\alpha_i}}$ is nontrivial
for $i=8,9,10$. From (\ref{eq:comrel8910}), interpreted in
characteristic~$2$, we see that 
$(\lambda_{0,0,0,0,d,e,f})^{x_1(r) x_2(s) x_4(t)}=\lambda_{0,0,0,0,d,e,f}$
implies 
\begin{eqnarray*}
\phi(dd_8+ed_9+fd_{10}) & = & \phi((drs+ert+fst)d_3\\
&& +(ds+et)d_5+(dr+ft)d_6\\
&& +(er+fs)d_7+dd_8+ed_9+fd_{10})
\end{eqnarray*}
for all $d_i \in \F_q$. Thus, if 
$(\lambda_{0,0,0,0,d,e,f})^{x_1(r) x_2(s) x_4(t)}=\lambda_{0,0,0,0,d,e,f}$
then the coefficients of $d_3, d_5, d_6, \dots, d_{10}$ imply
$r=s=t=0$ and so $\Stab_{\overline{K}}(\lambda_{0,0,0,0,d,e,f})=\{1\}$.

Now, suppose 
$(\lambda_{x,0,0,c,d,e,f})^{x_1(r) x_2(s) x_4(t)}=\lambda_{x',0,0,c',d',e',f'}$ 
where $x \in \F_q$ and $c,c',d,d',e,e',f,f' \in \F_q^\times$. 
From (\ref{eq:comrel8910}) we see that 
\[
(\lambda_{x,0,0,c,d,e,f})^{x_1(r) x_2(s)
  x_4(t)}=\lambda_{x',0,0,c',d',e',f'}
\]
if and only if
\begin{eqnarray*}
\phi(x'd_3+c'd_7+d'd_8+e'd_9+f'd_{10}) & = & \phi((x+ct+drs+ert+fst)d_3\\
&& \hspace{-0.25cm} +(ds+et)d_5+(dr+ft)d_6\\
&& \hspace{-0.25cm} +(c+er+fs)d_7+dd_8+ed_9+fd_{10})
\end{eqnarray*}
for all $d_i \in \F_q$. So 
$(\lambda_{x,0,0,c,d,e,f})^{x_1(r) x_2(s) x_4(t)}=\lambda_{x',0,0,c',d',e',f'}$
if and only if $d=d'$, $e=e'$, $f=f'$, $c=c'$, $r=\frac{f}{d}t$, 
$s=\frac{e}{d}t$ and $x'=x+ct+\frac{ef}{d}t^2$. It follows that 
(\ref{eq:repseteven}) is indeed a set of representatives and
$|\Stab_{\overline{K}}(\lambda_{x,0,0,c,d,e,f})|=2$.
By Clifford theory, the induced characters 
$\psi_{a_1,a_2,a_3} := \lambda_{0,0,0,0,a_1,a_2,a_3}^{\overline{U}}$
are $(q-1)^3$ distinct irreducible characters of degree $q^3$ 
of $\overline{U}$.
Furthermore, each $\lambda_{x,0,0,c,d,e,f}$ extends in two ways to its
inertia subgroup and by inducing to $\overline{U}$ we obtain
$4(q-1)^4$ irreducible characters $\psi_{x,a_1,a_2,a_3,a_4}$ of degree
$\frac{q^3}{2}$ of $\overline{U}$. By inflation, we obtain 
$\chi_{8,9,10,q^3,a_1,a_2,a_3} \in \Irr(U)$ of degree $q^3$ of~$U$ 
and $\chi_{8,9,10,\frac{q^3}{2},x,a_1,a_2,a_3,a_4}$ of degree
$\frac{q^3}{2}$ of $U$.

\medskip

\noindent \textbf{The irreducible characters in families
  $\mathcal{F}_{8,9}$, $\mathcal{F}_{8,10}$, $\mathcal{F}_{9,10}$.} 
The characters\linebreak in family $\mathcal{F}_{8,9}$ are those 
$\chi \in \Irr(U)$ such that 
$X_{\alpha_i} \subseteq \ker(\chi)$ for $i=10,11,12$ 
and $X_{\alpha_j} \not \subseteq \ker(\chi)$ for $j=8,9$.
We are going to work in the factor group 
$\overline{U} := U/X_{\alpha_{10}}X_{\alpha_{11}}X_{\alpha_{12}}$. 
By the commutator relations, the group
\[
\overline{N}_{10} := H_{\alpha_8}\prod_{i \ge 9} X_{\alpha_i}/X_{\alpha_{10}}X_{\alpha_{11}}X_{\alpha_{12}} 
\]
is a normal subgroup of $\overline{U}$ and we have
\[
\overline{N}_{10} \cong H_{\alpha_8} \times X_{\alpha_9} \quad
\text{and} \quad \overline{U}/\overline{N}_{10} \cong H_{\alpha_7},
\]
where $H_{\alpha_7}$, $H_{\alpha_8}$ are special $p$-groups of type
$q^{1+2}$ and $q^{1+4}$, respectively. So $\overline{N}_{10}$ has
$(q-1)^2$ almost faithful irreducible characters and they have degree
$q^2$. We claim that all of them extend to $\overline{U}$. 
Consider the subgroup 
\[
\overline{K}_{10} := \prod_{\genfrac{}{}{0pt}{}{i=1}{i \neq 1,5}}^{12}
X_{\alpha_i}/X_{\alpha_{10}}X_{\alpha_{11}}X_{\alpha_{12}}
\]
of $\overline{U}$ of index $q^2$ and let $a_1, a_2 \in \F_q^\times$.
Since $\overline{K}_{10} \cong
X_{\alpha_2}X_{\alpha_3}X_{\alpha_4}X_{\alpha_6}X_{\alpha_7} \times
X_{\alpha_8} \times X_{\alpha_9}$ there is a linear character
$\lambda_{a_1,a_2}$ of $\overline{K}_{10}$ such that 
\[
\lambda_{a_1,a_2}|_{X_{\alpha_8}}=\varphi_{\alpha_8,a_1} \quad \text{and}
\quad \lambda_{a_1,a_2}|_{X_{\alpha_9}}=\varphi_{\alpha_9,a_2}
\]
(here and in the following we identify the root subgroups with
their images in~$\overline{U}$). The induced characters 
$\lambda_{a_1,a_2}^{\overline{U}}$ have degree $q^2$ and are irreducible,
because they restrict irreducibly to $H_{\alpha_8}$. Hence they are
extensions of the almost faithful irreducible characters of
$\overline{N}_{10}$. Applying Gallagher's theorem to
$\overline{N}_{10} \unlhd \overline{U}$, we obtain 
$(q-1)^3$ almost faithful irreducible characters
$\psi_{8,9,q^3,a_1,a_2,a_3}$ of degree $q^3$ and $q^2(q-1)^2$ almost
faithful irreducible characters $\psi_{8,9,q^2,a_1,a_2,b_1,b_2}$ of
degree $q^2$ of $\overline{U}$, where $a_i \in \F_q^\times$ and 
$b_i \in \F_q$. By inflation, we obtain the irreducible characters 
$\chi_{8,9,q^3,a_1,a_2,a_3}$ and $\chi_{8,9,q^2,a_1,a_2,b_1,b_2}$ of $U$.
This completes the construction of the irreducible characters in the
family $\mathcal{F}_{8,9}$. 

The characters in the family $\mathcal{F}_{8,10}$ are those 
$\chi \in \Irr(U)$ such that 
$X_{\alpha_i} \subseteq \ker(\chi)$ for $i=9,11,12$ 
and $X_{\alpha_j} \not \subseteq \ker(\chi)$ for $j=8,10$,
and the characters in the family $\mathcal{F}_{9,10}$ are those 
$\chi \in \Irr(U)$ such that 
$X_{\alpha_i} \subseteq \ker(\chi)$ for $i=8,11,12$ 
and $X_{\alpha_j} \not \subseteq \ker(\chi)$ for $j=9,10$.
The definition and construction of the irreducible characters in these
two families are analogous to the definition and construction of the
irreducible characters in $\mathcal{F}_{8,9}$.

\medskip

\noindent \textbf{The irreducible characters in families
  $\mathcal{F}_{8}$, $\mathcal{F}_{9}$, $\mathcal{F}_{10}$.} 
The characters in family $\mathcal{F}_8$ are those 
$\chi \in \Irr(U)$ such that 
$X_{\alpha_i} \subseteq \ker(\chi)$ for $i=9,10,11,12$ 
and $X_{\alpha_8} \not \subseteq \ker(\chi)$.
We are going to work in the factor group 
$\overline{U} :=
U/X_{\alpha_9}X_{\alpha_{10}}X_{\alpha_{11}}X_{\alpha_{12}}$. 
By the commutator relations, the group
\[
\overline{N}_{9,10} := H_{\alpha_8}\prod_{i \ge 9}
X_{\alpha_i}/X_{\alpha_9}X_{\alpha_{10}}X_{\alpha_{11}}X_{\alpha_{12}}  
\]
is a normal subgroup of $\overline{U}$ and 
$\overline{U}/\overline{N}_{9,10} \cong H_{\alpha_7}$.
Since $\overline{N}_{9,10}$ is special of type $q^{1+4}$, it has $q-1$
almost faithful irreducible characters and they have degree $q^2$ and
the midafis $\mu_{\alpha_8,a}$, $a \in \F_q^\times$, are extensions of
these almost faithful characters to $\overline{U}$.
Now, Gallagher's theorem applied to 
$\overline{N}_{9,10} \unlhd \overline{U}$ gives us $(q-1)^2$
almost faithful irreducible characters $\psi_{8,q^3,a_1,a_2}$ of
degree $q^3$ and $q^2(q-1)$ almost faithful irreducible characters
$\psi_{8,q^2,a,b_1,b_2}$ of degree $q^2$ of $\overline{U}$, where
$a, a_i \in \F_q^\times$ and $b_i \in \F_q$. By inflation, we obtain the
irreducible characters $\chi_{8,q^3,a_1,a_2}$ and
$\chi_{8,q^2,a,b_1,b_2}$ of $U$. This completes the construction of
the irreducible characters in the family $\mathcal{F}_{8}$. 

The characters in the family $\mathcal{F}_{9}$ are those 
$\chi \in \Irr(U)$ such that 
$X_{\alpha_i} \subseteq \ker(\chi)$ for $i=8,10,11,12$ 
and $X_{\alpha_9} \not \subseteq \ker(\chi)$,
and the characters in the family $\mathcal{F}_{10}$ are those 
$\chi \in \Irr(U)$ such that 
$X_{\alpha_i} \subseteq \ker(\chi)$ for $i=8,9,11,12$ 
and $X_{\alpha_{10}} \not \subseteq \ker(\chi)$.
The definition and construction of the irreducible characters in these
two families are analogous to the definition and construction of the
irreducible characters in $\mathcal{F}_8$.

\medskip

\noindent \textbf{The irreducible characters in family $\mathcal{F}_{5,6,7}$.} 
The characters in $\mathcal{F}_{5,6,7}$ are those 
$\chi \in \Irr(U)$ such that 
$X_{\alpha_i} \subseteq \ker(\chi)$ for $i=8,9,\dots,12$ 
and $X_{\alpha_j} \not \subseteq \ker(\chi)$ for $j=5,6,7$.
We are going to work in the factor group 
$\overline{U} := U/\prod_{i \ge 8} X_{\alpha_i}$. 
By the commutator relations, the group
\[
\overline{N}_{5,6,7} := \prod_{\genfrac{}{}{0pt}{}{i=1}{i \neq 2,4}}^{12} X_{\alpha_i}/\prod_{i \ge 8} X_{\alpha_i}
\]
is a normal subgroup of $\overline{U}$ and we have
\[
\overline{N}_{5,6,7} \cong H_{\alpha_5} \times X_{\alpha_6} \times X_{\alpha_7} \quad
\text{and} \quad \overline{U}/\overline{N}_{5,6,7} \cong X_{\alpha_2}
\times X_{\alpha_4},
\]
where $H_{\alpha_5}$ is a special $p$-group of type $q^{1+2}$. 
Hence there are $(q-1)^3$ almost faithful irreducible characters of
$\overline{N}_{5,6,7}$ such that $X_{\alpha_5}$, $X_{\alpha_6}$,
$X_{\alpha_7}$ are not contained in their kernel, these characters
have degree $q$ and can be parametrized by $\psi_{a_1,a_2,a_3}$ where
$a_i \in \F_q^\times$. 

We show that each $\psi_{a_1,a_2,a_3}$ extends to
$\overline{U}$: Since 
$X_{\alpha_1} X_{\alpha_2} \cdot \prod_{i \ge 4} X_{\alpha_i} / \prod_{i \ge 8} X_{\alpha_i}$
is a normal abelian subgroup of $\overline{U}$ of index $q$, 
we can conclude from Ito's theorem \cite[Theorem (6.15)]{Isaacs:76}
that the degrees of all irreducible characters of $\overline{U}$
divide $q$. So from Clifford theory and Gallagher's theorem applied
to $\overline{N}_{5,6,7} \unlhd \overline{U}$ it follows that each 
$\psi_{a_1,a_2,a_3}$ extends to $\overline{U}$ in $q^2$ ways. We
denote these extensions by $\psi_{a_1,a_2,a_3,b_1,b_2}$ where 
$a_i \in \F_q^\times$ and $b_j \in \F_q$. Inflating 
$\psi_{a_1,a_2,a_3,b_1,b_2}$ to $U$, we get the $q^2(q-1)^3$
irreducible characters $\chi_{a_1,a_2,a_3,b_1,b_2}$ of degree $q$ of~$U$. 

\medskip

\noindent \textbf{The irreducible characters in family $\mathcal{F}_{5,6}$,
  $\mathcal{F}_{5,7}$, $\mathcal{F}_{6,7}$.}  
The characters in $\mathcal{F}_{5,6}$ are those 
$\chi \in \Irr(U)$ such that 
$X_{\alpha_i} \subseteq \ker(\chi)$ for $i \ge 7$ 
and $X_{\alpha_j} \not \subseteq \ker(\chi)$ for $j=5,6$.
Let $\overline{U}$, $\overline{N}_{5,6,7}$ be the groups defined in
the construction of the irreducible characters in the family
$\mathcal{F}_{5,6,7}$. The group $\overline{N}_{5,6,7} \cong
H_{\alpha_5} \times X_{\alpha_6} \times X_{\alpha_7}$ has $(q-1)^2$
irreducible characters such that $X_{\alpha_7}$ is contained in their
kernel and $X_{\alpha_5}$, $X_{\alpha_6}$ are not contained in their kernel,
these characters have degree $q$ and can be parametrized by
$\psi_{a_1,a_2}$ where $a_1, a_2 \in \F_q^\times$. 
Using Gallagher's theorem in the same way as for $\mathcal{F}_{5,6,7}$
we see that each $\psi_{a_1,a_2}$ extends to $\overline{U}$ in $q^2$
ways leading to the irreducible characters 
$\chi_{a_1,a_2,b_1,b_2}$ of degree~$q$ of~$U$.  

The characters in the family $\mathcal{F}_{5,7}$ are those 
$\chi \in \Irr(U)$ such that 
$X_{\alpha_i} \subseteq \ker(\chi)$ for $i=6,8,9,10,11,12$ 
and $X_{\alpha_j} \not \subseteq \ker(\chi)$ for $j=5,7$,
and the characters in the family $\mathcal{F}_{6,7}$ are those 
$\chi \in \Irr(U)$ such that 
$X_{\alpha_i} \subseteq \ker(\chi)$ for $i=5,8,9,10,11,12$ 
and $X_{\alpha_j} \not \subseteq \ker(\chi)$ for $j=6,7$.
The definition and construction of the irreducible characters in these
two families are analogous to the definition and construction of the
irreducible characters in $\mathcal{F}_{5,6}$.

\medskip

\noindent \textbf{The irreducible characters in families
  $\mathcal{F}_{5}$, $\mathcal{F}_{6}$, $\mathcal{F}_{7}$.} 
The characters in $\mathcal{F}_5$ are those 
$\chi \in \Irr(U)$ such that 
$X_{\alpha_i} \subseteq \ker(\chi)$ for $i \ge 6$ 
and $X_{\alpha_5} \not \subseteq \ker(\chi)$.
Let $\overline{U}$, $\overline{N}_{5,6,7}$ be the groups defined in
the construction of the irreducible characters in the family
$\mathcal{F}_{5,6,7}$. The group
$\overline{N}_{5,6,7} \cong H_{\alpha_5} \times X_{\alpha_6} \times
X_{\alpha_7}$ has $q-1$ irreducible characters 
such that $X_{\alpha_6}$ and $X_{\alpha_7}$ are contained in their
kernel and $X_{\alpha_5}$ is not contained in their kernel,
these characters have degree $q$ and can be parametrized by
$\psi_a$ where $a \in \F_q^\times$. 
Using Gallagher's theorem in the same way as for $\mathcal{F}_{5,6,7}$
we see that each $\psi_{a_1,a_2}$ extends to $\overline{U}$ in $q^2$
ways leading to the irreducible characters 
$\chi_{a,b_1,b_2}$ of degree~$q$ of~$U$.  

The characters in the family $\mathcal{F}_6$ are those 
$\chi \in \Irr(U)$ such that 
$X_{\alpha_i} \subseteq \ker(\chi)$ for $i=5,7,8,9,10,11,12$ 
and $X_{\alpha_6} \not \subseteq \ker(\chi)$,
and the characters in the family $\mathcal{F}_7$ are those 
$\chi \in \Irr(U)$ such that 
$X_{\alpha_i} \subseteq \ker(\chi)$ for $i=5,6,8,9,10,11,12$ 
and $X_{\alpha_7} \not \subseteq \ker(\chi)$.
The definition and construction of the irreducible characters in these
two families are analogous to the definition and construction of the
irreducible characters in $\mathcal{F}_5$.

\medskip

\noindent \textbf{The irreducible characters in family $\mathcal{F}_{lin}$.} 
We have 
\[
U/[U,U] \cong X_{\alpha_1} \times X_{\alpha_2} \times
X_{\alpha_3} \times X_{\alpha_4} \cong \F_q^4
\]
and so the statements about the linear characters of $U$ are
clear. This completes the proof of Theorem~\ref{thm:irrU}. 
\hfill $\Box$ 

\subsection{Numbers and degrees of irreducible characters}
\label{subsec:nrsdegs}

We see from Table~\ref{tab:irrU} that for odd $q$ the degrees of the
irreducible characters of $U$ are powers of $q$. So~$U$ is a
$q$-power-degree group in the sense of~\cite{IsaacsCharAlg}. This
observation is a special case of a general theorem of B.~Szegedy on
the Sylow $p$-subgroups of classical groups defined over finite fields
of of good characteristic $p$, see~\cite[Theorem~2]{Szegedy}. 

The subgroup $U_n(q)$ of $\GL_n(\F_q)$ consisting of all upper
unitriangular matrices is a Sylow $p$-subgroup of $\GL_n(\F_q)$. A
conjecture of G.~Higman~\cite{Higman} states that the number of
conjagacy classes of $U_n(q)$ is given by a polynomial in $q$ with
integer coefficients. I.M.~Isaacs proved that the degrees of the
irreducible characters of $U_n(q)$ are of the form $q^e$, 
$0 \le e \le \mu(n)$ where the upper bound $\mu(n)$ depends on $n$ and
is known explicitly. G.~Lehrer~\cite{Lehrer} conjectured that the
number $N_{n,e}(q)$ of irreducible characters of $U_n(q)$ of degree
$q^e$ are given by a polynomial in $q$ with integer
coefficients. I.M.~Isaacs suggested a strengthened form of
Lehrer's conjecture stating that $N_{n,e}(q)$ is given by a polynomial
in $q-1$ with nonnegative integer coefficients. Obviously, Isaac's
conjecture implies Higman's and Lehrer's conjectures.

From Table~\ref{tab:irrU}, we can derive that an analogue of Isaac's
conjecture holds for the Sylow $p$-subgroup $U$ of the Chevalley
groups of type $D_4$, even in bad characteristic.

\begin{corollary} \label{cor:degs}
The degrees of the irreducible characters of $U$ are given by
Table~\ref{tab:irrUdegs}. In particular, the number of conjugacy
classes of $U$ is 
\[
\begin{cases}2q^5+5q^4-4q^3-4q^2+2q & \text{, if   } $q$ \text{   is odd,}\\
2q^5+8q^4-16q^3+14q^2-10q+3 & \text{, if   } $q$ \text{   is even.}
\end{cases}
\]
\begin{table}[!ht] 
\caption{Numbers and degrees of the irreducible characters of $U$. The
  numbers of the irreducible characters are given as polynomials in $v=q-1$.} 
\label{tab:irrUdegs}

\begin{center}
\begin{tabular}{l|l|l}
\hline
\rule{0cm}{0.4cm}
Degree & Number of irreducible characters & Comments
\rule[-0.1cm]{0cm}{0.4cm}\\
\hline
\rule{0cm}{0.5cm}
$q^4$ & $v^4+3v^3+3v^2+v$ &
\rule[-0.2cm]{0cm}{0.4cm}\\
\hline
\rule{0cm}{0.4cm}
$q^3$ & $v^5+5v^4+10v^3+7v^2+v$ & if $q$ is odd\\
      & $v^5+4v^4+10v^3+7v^2+v$ & if $q$ is even
\rule[-0.2cm]{0cm}{0.4cm}\\
\hline
\rule{0cm}{0.4cm}
$\frac{q^3}{2}$ & $4v^4$ & only if $q$ is even
\rule[-0.2cm]{0cm}{0.4cm}\\
\hline
\rule{0cm}{0.4cm}
$q^2$ & $3v^4+9v^3+9v^2+3v$ &
\rule[-0.2cm]{0cm}{0.4cm}\\
\hline
\rule{0cm}{0.4cm}
$q$ & $v^5+5v^4+10v^3+9v^2+3v$ &
\rule[-0.2cm]{0cm}{0.4cm}\\
\hline
\rule{0cm}{0.4cm}
$1$ & $v^4+4v^3+6v^2+4v+1$ &
\rule[-0.2cm]{0cm}{0.4cm}\\
\hline
\end{tabular}
\end{center}
\end{table}
\end{corollary}
\begin{proof}
This follows from Theorem~\ref{thm:irrU}.
\end{proof}
For odd prime powers $q$ the number of conjugacy classes of $U$ was already
computed by S.M.~Goodwin and G.~R\"ohrle~\cite[Table~1]{GoodwinRoehrle}.

\bigskip

\textbf{Acknowledgements.} Part of this work was done while the authors
were participating in the program on Representation Theory of Finite
Groups and Related Topics at the Mathematical Sciences Research
Institute (MSRI), Berkeley. It is a pleasure to thank the organizers
Professors J. L.~Alperin, M.~Brou\'e, J. F.~Carlson, A. S.~Kleshchev,
J.~Rickard, B.~Srinivasan for generous hospitality and support.
We also thank C.~Andr{\'e}, P.~Diaconis, I.M.~Isaacs, N.~Thiem and
N.~Yan for the stimulating seminar on supercharacters at the MSRI.

%%%%%%%%%%%%%%%%%%%%%%%%%%%%%%%%%%%%%%%%%%%%%%%%%%%%%%%%%%%%%%%%%%%%%%%%%%%%%%%%

%\bibliographystyle{mystyle}
%\bibliography{mybib}

\end{document}